# Analytical Identification of Design and Multidimensional Spaces Using R-Functions

S. Kucherenko[a*], O. Klymenko[b], N. Shah[a]
[a]Imperial College London, London, SWT 2AZ, UK
s.kucherenko@imperial.ac.uk
[b]University of Surrey, Guildford GU2 7XH, UK

**Abstract.** The design space (DS) is defined as the combination of materials and process conditions that guarantees the assurance of quality. This principle ensures that as long as a process operates within DS, it consistently produces a product that meets specifications. A novel DS identification method called the R-DS identifier has been developed. It makes no assumptions about the underlying model - the only requirement is that each model constraint (CQA) should be approximated by a multivariate polynomial model. The method utilizes the methodology of Rvachev's R-functions and allows for explicit analytical representation of the DS with only a limited number of model runs. R-functions provide a framework for representing complex geometric shapes and performing operations on them through implicit functions. The theory of R-functions enables the solution of geometric problem such as identification of DS through algebraic manipulation. It is more practical than traditional sampling or optimization-based methods. The method is illustrated using a batch reactor system.

**Key words:** Design Space, R-functions, Multivariate Polynomial Model

## 1. Introduction

René Descartes introduced the method of coordinates to connect geometric objects (points, lines, bodies) with analytical objects (numbers, equations, inequalities). Historically, significant focus has been placed on the direct problem of analytic geometry: the study of curves and surfaces defined by equations.

The inverse problem, that is expressing geometric objects with equations, has been solved for basic shapes like straight lines, circles, conic sections, and second-degree surfaces, which are usually described by polynomial equations. Although polynomials are relatively simple, solving these equations can be very challenging, seemingly restricting the inverse problem's application to more complex objects. However, by incorporating additional operations, the analytical methods' descriptive power can be greatly enhanced. This enhancement allows for the creation of equations for geometric objects of nearly any shape, thereby solving the inverse problem of analytic geometry in principle.

Consider two circles described by the equations:
$$R_0^2 = (x - x_0)^2 + (y - y_0)^2,$$
$$R_1^2 = (x - x_1)^2 + (y - y_1)^2. \qquad (1)$$

They are shown for illustration in Fig. 1. We aim to analytically describe the overlapping region of the areas inside the circles, given by $\phi_1(x, y) \geq 0, \phi_2(x, y) \geq 0$, where implicit functions $\phi_1(x, y), \phi_2(x, y)$ have the form
$$\phi_1(x, y) = R_0^2 - (x - x_0)^2 - (y - y_0)^2,$$
$$\phi_2(x, y) = R_1^2 - (x - x_1)^2 - (y - y_1)^2. \qquad (2)$$

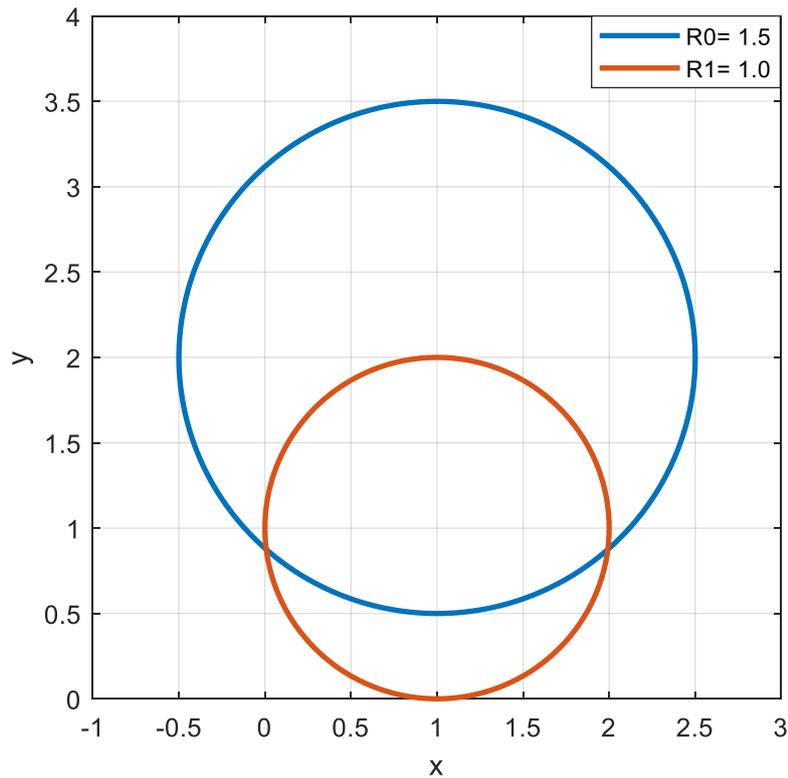

Fig. 1. Two circles. Test case 4.1.

Is it possible to derive analytical formulas for the boundary of this overlapping area given the implicit functions (2)? We recall that an implicit surface representation has the form $f(x) = 0$, where $x$ is a point on the surface implicitly described by the function $f(x)$. Can we perform operations on functions describing implicit surfaces using logical operations such as conjunction (AND), disjunction (OR), and NOT? (Fig. 2)? The answer is "yes", it can be done using R-functions.

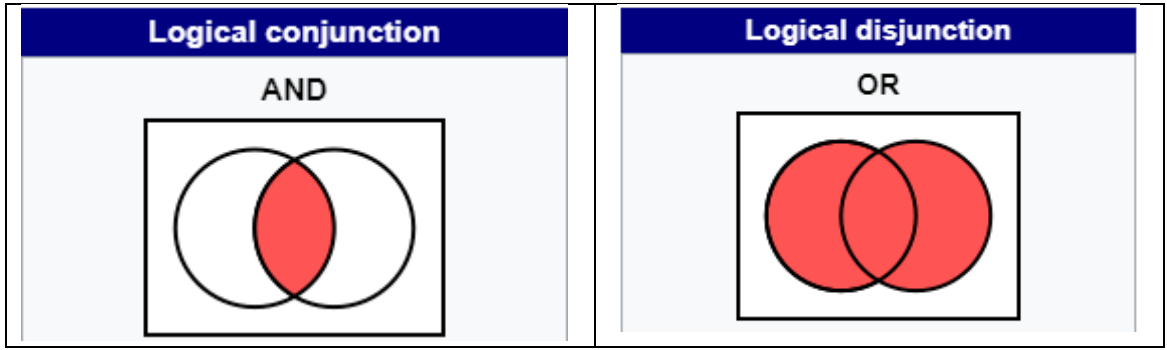

Fig. 2. Boolean functions.

## 2. Introduction to R-Functions

Rvachev's R-functions, named after the Ukrainian mathematician Vladimir Rvachev (Rvachev, 1982, Shapiro, 1991) are a powerful mathematical tool used in various fields, including computer graphics, computational geometry, and solid modelling. R-functions provide a formal framework for representing complex geometric shapes and performing operations on them through implicit functions. This introduction explores the fundamental concepts, applications, and advantages of R-functions.

R-functions are built on the idea of using implicit functions to define geometric objects. An implicit function $f(x)$ represents a shape such that the points inside the shape satisfy $f(x) \geq 0$, points on the boundary satisfy $f(x) = 0$, and points outside satisfy $f(x) \leq 0$. R-functions extend this concept by enabling the combination of these implicit functions to perform set-theoretic operations like union, intersection, and difference.

### 2.1. Basic R-Functions

The basic R-functions are defined using logical operations on implicit functions:

- Union (OR operation): The union of two implicit functions $f_1$ and $f_2$ is given by an R-function that satisfies $f_\vee = f_{\text{union}}(x) \geq 0$ if either $f_1(x) \geq 0$ or $f_2(x) \geq 0$.
- Intersection (AND operation): The intersection is defined by an R-function that satisfies $f_\wedge = f_{\text{intersection}}(x) \geq 0$ if both $f_1(x) \geq 0$ and $f_2(x) \geq 0$.
- Difference (NOT operation): The difference between two implicit functions is represented by an R-function that modifies $f_1$ based on the complement of $f_2$.

These operations are achieved through specific algebraic combinations of the implicit functions, ensuring the resulting function maintains the desired geometric properties.

### 2.2. Applications of R-Functions

R-functions have a wide range of applications across various domains:

- **Computer Graphics**: In computer graphics, R-functions are used for procedural modelling, enabling the creation of complex shapes and structures through simple algebraic expressions.
- **Computational Geometry**: R-functions facilitate the representation and manipulation of geometric objects, making them useful in algorithms for collision detection, mesh generation, and object recognition.
- **Solid Modelling**: In CAD systems, R-functions provide a robust method for constructing and editing solid models, supporting operations like Boolean combinations and freeform shape design.
- **Engineering Analysis**: R-functions are employed in finite element analysis and other engineering simulations to define domains and boundary conditions with high precision.

## 2.3. Advantages of R-Functions

The use of R-functions offers several advantages:

- **Compact Representation**: Complex shapes can be represented using relatively simple algebraic expressions, reducing storage and computational requirements.
- **Smooth Transitions**: R-functions inherently provide smooth transitions between combined shapes, which is beneficial for applications requiring high-quality surface continuity.
- **Flexibility**: The framework allows for easy modification and combination of shapes, making it highly versatile for various design and analysis tasks.
- **Mathematical Rigour**: R-functions offer a mathematically rigorous way to perform set-theoretic operations on geometric objects, ensuring consistency and accuracy.

## 3. Theory of R-Functions

Consider a function $S(x)$ defined on the real axis as follows (further we follow notations of Shapiro, 1991):

$$S(x) = \begin{cases} 0 & \text{if } x \leq -0, \\ 1 & \text{if } x \geq +0. \end{cases} \quad (3)$$

A function $y = f(x_1, \ldots, x_n)$ is an $R$-function if there exists a Boolean function $Y = F(X_1, \ldots, X_n)$ such that the following equality is satisfied:

$$S[(f(x_1, x_2, \ldots, x_n)] = F[S(x_1), S(x_2), \ldots, S(x_n)]. \quad (4)$$

We will refer to this Boolean function as a companion function of a given $R$-function. Informally, a real function is an $R$-function if it can change its property (sign) only when some of its arguments change the same property (sign).

The following functions are $R$-functions with their corresponding Boolean companion function in parentheses:

$$\begin{aligned} C &\equiv \text{const}; & \text{(logical 1)} \\ \bar{x} &\equiv -x; & \text{(logical negation } \neg \text{)} \\ x_1 \wedge_1 x_2 &\equiv \min(x_1, x_2); & \text{(logical conjunction } \wedge \text{)} \\ x_1 \vee_1 x_2 &\equiv \max(x_1, x_2); & \text{(logical disjunction } \vee \text{)} \end{aligned} \quad (5)$$

The set of $R$-functions is infinite. However, for applications, it is not necessary to know all $R$-functions. This leads to the notion of sufficiently complete systems of $R$-functions.

**Theorem 1.** Let $H$ be some system of $R$-functions, and $G$ be the corresponding system of companion Boolean functions. The system $H$ is sufficiently complete, if the system $G$ is complete.

For example, take $G = \{0, \neg X, X_1 \wedge X_2\}$. It is well known that all logic functions can be constructed using just conjunction and negation; in other words, $G$ is complete. $R$-negation and $R$-conjunction are defined as

$$\begin{aligned} \bar{x} &\equiv -x; \\ x_1 \wedge_0 x_2 &\equiv x_1 + x_2 - \sqrt{x_1^2 + x_2^2}. \end{aligned} \quad (6)$$

For practical reasons it is useful to define a more general set $R_\alpha$ of R-functions:

**Definition**: A system $R_\alpha$ is defined as follows:

$$\begin{aligned} x_1 \wedge_\alpha x_2 &\equiv \frac{1}{1+\alpha}\left(x_1 + x_2 - \sqrt{x_1^2 + x_2^2 - 2\alpha x_1 x_2}\right); \\ x_1 \vee_\alpha x_2 &\equiv \frac{1}{1+\alpha}\left(x_1 + x_2 + \sqrt{x_1^2 + x_2^2 - 2\alpha x_1 x_2}\right), \end{aligned} \quad (7)$$

where $\alpha(x_1, x_2)$ is an arbitrary symmetric function such that $-1 < \alpha(x_1, x_2) \leq 1$.

The first function is called **R**-conjunction, while the second function - **R**-disjunction. Further we will be using the case of $\alpha = 1$ and a set of $x, y$ variables. Then (7) will have the form:

$$\begin{aligned} x \wedge_1 y &= \frac{1}{2}\left(x + y - (x^2 + y^2 - 2xy)^{\frac{1}{2}}\right), \\ x \vee_1 y &= \frac{1}{2}\left(x + y + (x^2 + y^2 - 2xy)^{\frac{1}{2}}\right). \end{aligned} \quad (8)$$

The set of all points in Euclidean space $E^n$ where a function $y = f(x_1, \ldots, x_n) = 0$ is called a drawing. The set of points where $f(x_1, \ldots, x_n) \geq 0$ is called a region. Many known geometric objects are algebraic drawings, making their corresponding polynomial functions "building blocks" or primitives. By combining these primitives with logical operations ("and" $\wedge$, "or" $\vee$, "not" $\neg$, etc.), most geometric objects of practical interest can be described.

## 3.1 R-functions: From Boolean expressions to real functions

Suppose we are given a description of a geometric object $D \subseteq E^n$

$$D = F[(\phi_1 \geq 0), \ldots, (\phi_m \geq 0)], \tag{9}$$

where real-function inequalities $\phi_i(x_1, \ldots, x_n) \geq 0$ define primitive geometric regions, and $F$ is a set function constructed using standard set operations ∩,∪,− on these primitive regions. Alternatively, replacing the set operations by the corresponding logical functions ∧,∨,¬, we can view $F$ as a Boolean logic function. The equation (9) becomes a predicate equation

$$F[S_2(\phi_1), \ldots, S_2(\phi_m)] = 1, \tag{10}$$

which holds for all points $p \in D \subseteq E^n$. We seek a single real-function inequality $f(x_1, \ldots, x_n) \geq 0$ that defines the composite object $D$.

**Theorem 2.** Let $F(X_1, \ldots, X_m)$ be a closing Boolean function that is the companion of a continuous $R$-function $f(x_1, \ldots, x_m)$. Then for any continuous real functions $\phi_i$, $i = 1, 2, \ldots, m$, and a closed region $D \subset E^n$ given by Eq. (9), $D$ can be defined by inequality

$$f(\phi_1, \ldots, \phi_m) \geq 0. \tag{11}$$

In other words, to obtain a real function defining the region $D$ constructed from primitive regions $\phi_i \geq 0$, it suffices to construct an appropriate $R$-function and substitute for its arguments the real functions $\phi_i$ defining the primitive regions.

## 4. Test cases. Geometrical objects

In this Section we consider two test cases in 2D. Test cases in 3D are considered in the Appendix.

### 4.1. Two circles

Consider two circles defined by Eq. (1), Fig. 1 and the corresponding implicit functions $\phi_1(x, y), \phi_2(x, y)$ defined by Eq. (2). We are interested in the intersection (AND / ∧ operation) and union (OR / ∨ operation) of the domains representing the interiors of these circles. Following the theory of R-functions, we would like to build continuous $R$-functions $f_\wedge(\phi_1, \phi_2) \geq 0, f_\vee(\phi_1, \phi_2) \geq 0$ from the set of continuous real functions $\phi_1(x, y), \phi_2(x, y)$. We follow the procedure outlined in the previous section: we build two objects $D_\wedge$ and $D_\vee$ defined by Eq. (11), where $f(\phi_1, \ldots, \phi_m)$ is defined by R-functions (8).

In the test we use the following parameters: $x_0 = 1$; $x_1 = 1$; $y_0 = 2$; $y_1 = 1$; $R_0 = 1.5$; $R_1 = 1$. Then

$$\begin{cases} f_\wedge(x, y) = x^2 - 3y - 2x + y^2 - \left(\frac{(8y-7)^2}{8}\right)^{\frac{1}{2}} + \frac{15}{8}, \\ f_\vee(x, y) = x^2 - 3y - 2x + y^2 + \left(\frac{(8y-7)^2}{8}\right)^{\frac{1}{2}} + \frac{15}{8}. \end{cases} \tag{12}$$

The boundaries of the intersection and union of the two circular domains are plotted in Fig. 3.

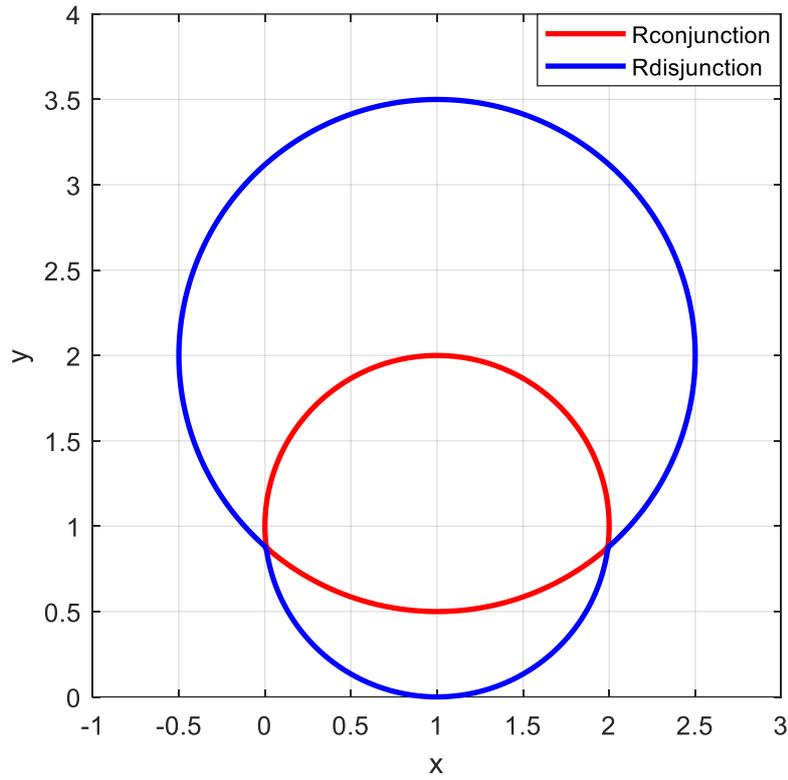

Fig. 3. Domain boundaries defined by implicit equations described by the $R$-functions (12): $f_\wedge(\phi_1, \phi_2) = 0$, $f_\vee(\phi_1, \phi_2) = 0$ for the test case 4.1.

### 4.2. Two parabolas

Consider two parabolas defined by implicit functions $\phi_1(x, y), \phi_2(x, y)$:
$$\phi_1(x, y) = y - a(x - x_0)^2 + d, \qquad (13)$$
$$\phi_2(x, y) = y + a(x - x_0)^2 + b.$$

They are shown in Fig. 4(a). We are interested in the intersection of two areas: $\phi_1(x, y) \geq 0, \phi_2(x, y) \leq 0$ and their union. They can be obtained as continuous $R$-functions $f_\wedge(\phi_1, \overline{\phi_2}) \geq 0, f_\vee(\phi_1, \overline{\phi_2}) \geq 0$:

$$\begin{cases} f_\wedge(x,y) = 2x - x^2 - \dfrac{((4y+9)^2)^{\left(\frac{1}{2}\right)}}{4} - \dfrac{1}{4}, \\ f_\vee(x,y) = 2x - x^2 + \dfrac{((4y+9)^2)^{\left(\frac{1}{2}\right)}}{4} - \dfrac{1}{4}. \end{cases} \qquad (14)$$

In the test we used the following parameters a= 1; $x_0$ = 1; $x_1$ = 1; $d$ = 3; b= 1.5.

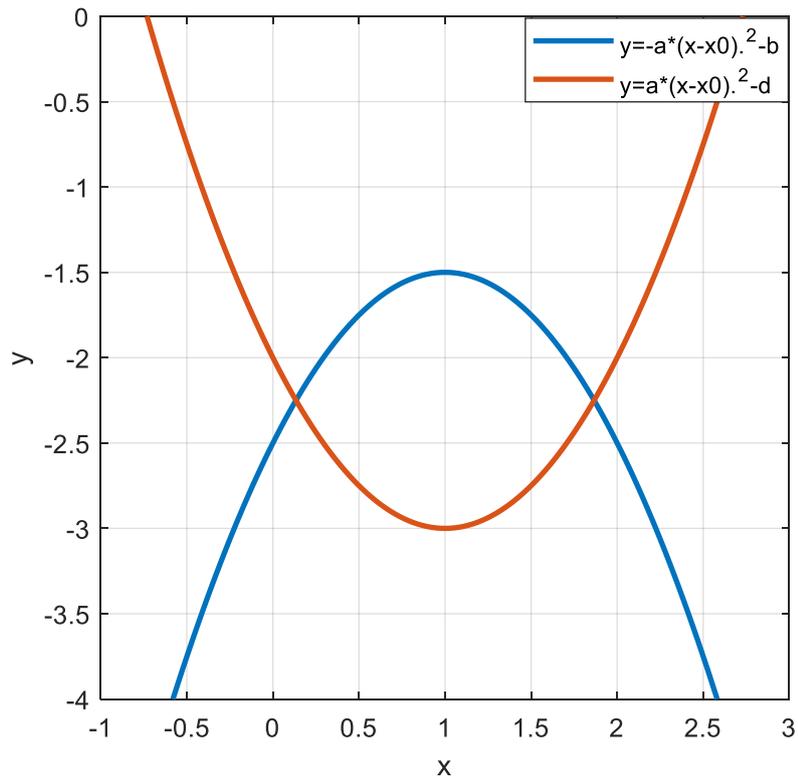

(a)

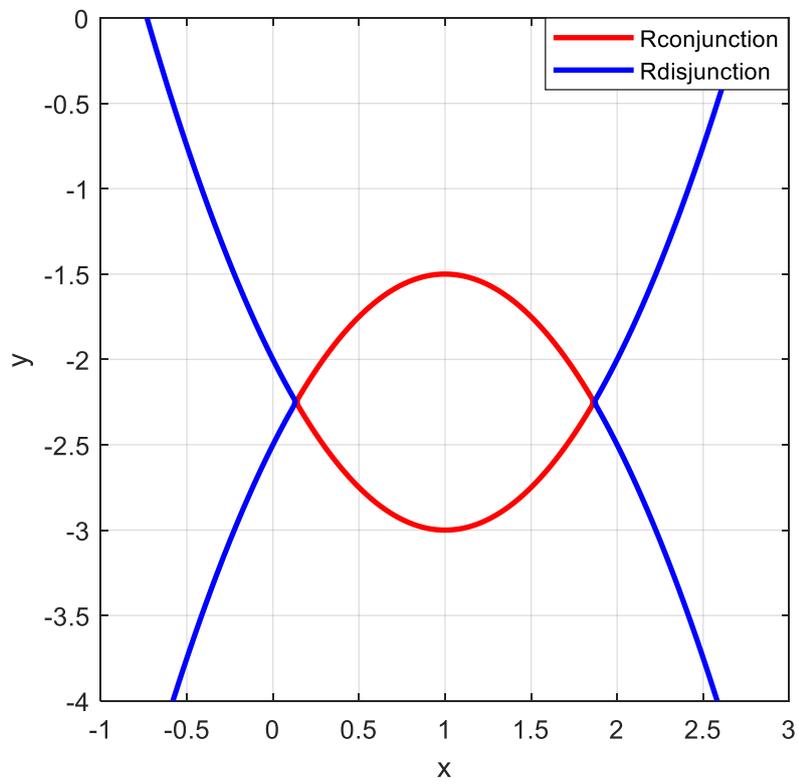

(b)

Fig. 4. Domain boundaries defined by implicit equations described by the $R$-functions (14): $f_\wedge(\phi_1, \phi_2) = 0$, $f_\vee(\phi_1, \phi_2) = 0$ for the test case 4.2.

## 5. Design space

The pharmaceutical manufacturing industry is subject to strict requirements to ensure drug quality and safety (Blacker and Williams, 2011). The modern approach to regulating pharmaceutical product development utilizes the concept of Quality by Design (QbD), introduced by the US Food and Drug Administration and the International Council for Harmonisation (ICH, 2009). QbD is based on the principle that quality should be built directly into the design of products and manufacturing processes. This approach ensures that the physical or chemical properties of a product, known as critical quality attributes (CQAs), remain within specified limits. CQAs can vary during processing due to inevitable variations in critical process parameters (CPPs), which are closely monitored and controlled. To manage all process variations, the QbD principle includes the use of the design space (DS), defined as "the multidimensional combination and interaction of input variables and process parameters that have been demonstrated to provide assurance of quality" (ICH, 2009). Within the DS, operators can adjust process parameter values without violating CQAs, providing manufacturers with operational flexibility.

Original sampling-based methods for defining the Design Space (DS) involved creating a dense grid of sample points across the process parameter space and performing simulations at each point to assess if CQAs are met, resulting in a deterministic DS. To calculate a probabilistic DS, a Monte Carlo simulation is conducted at each sampling point in the space of uncertain parameters. The probability is estimated as the fraction of times the Monte Carlo simulations produce results that comply with all CQAs. Although computationally demanding, this method provides a straightforward way to identify the DS (García-Muñoz et al., 2015).

Identification of probabilistic DS is a demanding task and for a typical practical problem the traditional approach based on exhaustive sampling requires costly computations. A novel theoretical and numerical framework for determining probabilistic DS using metamodelling and adaptive sampling was proposed in Kucherenko et al. (2020). It is based on the multi-step adaptive technique using a metamodel for a probability map as an acceptance-rejection criterion to optimize sampling for identification of the DS.

A challenge remains in the visualization of DS needed to communicate high dimensional probability maps and the inclusion of this information into the formal documentation required for the approval of DS by a government agency. For more than three process parameters visualization of DS is a challenging task. Guidance for Industry, Q8(R2) Pharmaceutical Development (2009) has the following recommendations: "When multiple parameters are involved, the design space can be presented for two parameters, at different values (e.g., high, middle, low) within the range of the third parameter, the fourth parameter, and so on." Considering DS in $d$-dimensional process parameters space, a graphical presentation of DS via 2D plots prescribed in Guidance for Industry, Q8(R2) Pharmaceutical Development (2009) would result in $N_P = \frac{d(d-1)}{2} 3^{d-2}$ 2D plots.

In all previous works, the DS was identified by considering all CQAs simultaneously. In this study, we take a different approach. Firstly, we identify the DS's separately for each constraint, aiming to find implicit functions to describe the DS boundaries analytically. In the second stage, we determine the joint DS as the intersection of all separate DS's using the theory of R-functions. This method yields the joint DS as an analytical expression, addressing also the task of communicating high-dimensional DS and including this information in formal documentation. Having the DS in analytical form allows for easy identification without running an underlying potentially expensive model and constraints whether a point in the process parameter space lies within the DS.

## 5.1. Formal definition of Design Spaces

Consider a model $f(\boldsymbol{u};\boldsymbol{\theta})$ and a vector of constraints $\mathbf{g}(\boldsymbol{u};\boldsymbol{\theta})$ defined as follows

$$y = f(\boldsymbol{u};\boldsymbol{\theta}), \tag{15}$$

$$\mathbf{g}(\boldsymbol{u};\boldsymbol{\theta}) \geq \mathbf{g}^*. \tag{16}$$

Here $\boldsymbol{u}$ is a vector of process parameters, $\boldsymbol{\theta}$ is a vector of uncertain parameters, $\mathbf{g}^*$ is a vector of the process thresholds (CQAs). The DS is defined as

$$DS = DS(\boldsymbol{u} \mid p(\boldsymbol{u}) \geq p^*), \tag{17}$$

$$p(\boldsymbol{u}) = E_\theta[I(\mathbf{g}(\boldsymbol{u},\boldsymbol{\theta}) \geq \mathbf{g}^*)]. \tag{18}$$

Eq. (17) and Eq. (18) define DS as a set of process parameters that ensures the probability of process meeting the constraints be greater than the critical probability threshold $p^*$. Here $I$ represents an indicator function taking value of either 1 (satisfying constraints) or 0 (not satisfying constraints).

Consider a case when a vector of uncertain model parameters $\boldsymbol{\theta}$ has a fixed value $\boldsymbol{\theta} = \boldsymbol{\theta}^*$. It is formally equivalent to a probability density function $\varphi(\boldsymbol{\theta})$ of $\boldsymbol{\theta}$ being equal to a delta function: $\varphi(\boldsymbol{\theta}) = \delta(\boldsymbol{\theta}^*)$. Then $p^*=1$ and DS is defined as

$$DS = DS(\boldsymbol{u}; p(\boldsymbol{u}) = 1), \tag{19}$$

or equivalently

$$DS = DS(\boldsymbol{u}; \mathbf{g}(\boldsymbol{u};\boldsymbol{\theta}^*) \geq \mathbf{g}^*). \tag{20}$$

This setting can also be seen as deterministic and as a subset of the general formulation (17),(18). In this work we consider only deterministic setting but the developed methodology can be easily extended for the case of probabilistic setting. Further we omit dependence of $\mathbf{g}(\mathbf{u}; \boldsymbol{\theta}^*)$ on the constant vector parameters $\boldsymbol{\theta}^*$ for brevity.

We highlight an important equality for the proposed methodology:

$$DS(\mathbf{u}; \mathbf{g}(\mathbf{u}) \geq \mathbf{g}^*) = DS\big(\mathbf{u}; \cap_{i=1}^{L} g_i(\mathbf{u}) \geq g_i^*\big) = \cap_{i=1}^{L} DS_i(\mathbf{u}; g_i(\mathbf{u}) \geq g_i^*). \quad (21)$$

## 5.2. The R-DS identifier method

To identify the DS the following algorithm is proposed:
1) construct Multivariate Polynomial Regression metamodels $g_i^M$, $i = 1, \ldots, K$ for each of the $K$ CQA's (constraints) $g_i(\mathbf{u})$, $i = 1, \ldots, K$ using $N$ model runs;
2) using a metamodel $g_i^M$ formulate $DS_i = DS_i(\mathbf{u}; g_i(\mathbf{u}) \geq g_i^*)$ via an implicit function $\phi_i(\mathbf{u})$: $DS_i = DS_i(\mathbf{u}; \phi_i(\mathbf{u}) \geq 0)$, $i = 1, \ldots, M$;
3) construct the joint $DS$:

$$DS = \cap_{i=1}^{L} DS_i \text{ as } DS = DS(\mathbf{u}, f_\wedge(\phi_i(\mathbf{u}), \ldots, \phi_M(\mathbf{u})) \geq 0). \quad (22)$$

This algorithm is much more efficient and practical than the existing sampling-based methods in that
1) metamodels are built for individual constraints and not for joint DS. In this case the number of required model runs $N$ is much smaller due to the smoothness of individual constraints, unlike complex shapes of the joint DS for which metamodels were built previously;
2) the identified joint DS is given in an explicit analytical form as an algebraic equation.

## 5.3. Batch Reactor System

We consider a two-stage batch reactor problem adapted from Samsatli et al. (1999), where the following chemical reaction takes place:

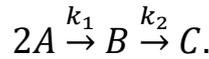

$$2A \xrightarrow{k_1} B \xrightarrow{k_2} C.$$

Here, B represents the desired product, while C denotes an undesired byproduct. The problem is formulated as the following system of DAE's with the initial conditions:

$$\begin{aligned}
\frac{dC_A'}{d\tau} &= t \cdot (-2 \cdot k_1 \cdot C_A{}^2), \\
\frac{dC_B}{d\tau} &= t \cdot (k_1 \cdot C_A{}^2 - k_2 \cdot C_B), \\
\frac{dC_C}{d\tau} &= t \cdot (k_2 \cdot C_B), \\
C_{A_1} &= C_A' + C_{A0}, \ C_{A0} = 2000, \\
k_j &= k_j^\circ e^{-E_j/RT}, \ j = 1,2.
\end{aligned} \quad (23)$$

$$C_A'(t = 0) = 0, \ C_B(t = 0),$$
$$C_C(t = 0) = 0.$$

The values of reaction kinetic parameters are given in Table 1.

Table 1. Values of the kinetic parameters.

| Parameter | | Value |
|---|---|---|
| $E_1$ (J/mol) | Activation energy | 2500.2 |
| $E_2$ (J/mol) | Activation energy | 5000.1 |
| $k_1^\circ$ (min$^{-1}$) | Pre-exponential factor | 0.0666 |
| $k_2^\circ$ (min$^{-1}$) | Pre-exponential factor | 10333.5 |

Two controlled variables (temperature $T$ and processing time $t$) can be adjusted independently to satisfy two CQA's: the product's *Purity* and *Profit* constraints: $Purity \geq P^* = 80\%$, $Profit \geq Pr^* = \$128/min$, where *Purity* and *Profit* are defined as

$$Purity = \frac{C_{B_1}(t)}{C_A(t) + C_B(t) + C_C(t)}, \quad (24)$$

$$Profit = \frac{(100 \cdot C_B(t) - 20 \cdot C_A) \cdot V}{t + 30}.$$

A reactor is a fixed volume vessel, with $V = 1$ m³. The ranges of controlled variables $\boldsymbol{u} = (t, T)$ are $T \in [250, 300]$ K, $t \in [250, 300]$ mins. The entire process operates isothermally.

Following the proposed methodology, we firstly need to define separately DS's for Purity and Profit

$$DS_{(\boldsymbol{u}; Purity(\vec{x}))} = (\boldsymbol{u}, Purity \geq P^*) \quad (25)$$

$$DS_{(\boldsymbol{u}; Profit(\boldsymbol{u}))} = (\boldsymbol{u}, Profit \geq Pr^*).$$

For this purpose we built Multivariate Polynomial Regression metamodels $g_i^M$, $i = 1, 2$ for each of the two constraints $g_i(x)$, $i = 1, 2$. *N* points are sampled in the space of controlled variables $\boldsymbol{u} = (t, T)$. For each sampled point the DAE system (23) is solved, and the input-output table is used to build metamodels.

In this work we used the SobolGSA software (Kucherenko and Zaccheus, 2025) to build multivariate polynomial models. There are two different methods implemented in SobolGSA: the Random Sampling-High Dimensional Model Representation (RS-HDMR) ( Li et al., 2002) and the Bayesian Sparce PCE method (BSPCE) (Shao et al., 2017). The RS-HDMR method is based on truncating the ANOVA-HDMR expansions up to the second or third order, while the truncated terms are then approximated by orthonormal polynomials (Zuniga et al. 2013). RS-HDMR belongs to a wider class of methods known as polynomial chaos expansion (PCE) ( Zuniga et al., 2013 ). PCE is a metamodeling method in which the output of a model is presented as a series of polynomial functions of the input variables, with each polynomial term weighted by coefficients that reflect the contribution of the corresponding input variable to the output variability. Usually only a few terms are really relevant in the PCE structure. The BSPCE method makes use of sparse PCE. Selection of the proposed PCE structure is based on a Bayesian approach using the Kashyap information criterion for model selection (Shao et al., 2017). It enables the construction of accurate metamodels using only a limited number of model runs. Sampling was performed using QMC sampling based on Sobol' sequences (Sobol et al., 2011). QMC sampling gives a much better way of arranging $N$ points in $d$–dimensions and a faster convergence rate than random sampling (Kucherenko et al., 2015).

Application of BSPCE revealed that for both constraints the best PCE approximation is a polynomial of the second order with respect to $T$ and the first order with respect to $t$:

$$g_i^M = a_{00}^{[i]} + a_{10}^{[i]}T + a_{20}^{[i]}T^2 + a_{01}^{[i]}t + a_{11}^{[i]}t\,T. \tag{26}$$

Coefficients of decomposition $\{a_{kl}^{[i]}\}$ are found by fitting a polynomial function (26) to a input-output data points. Tests showed that $N$ =64 model runs were sufficient to build metamodels for $Purity(t,T)$ and $Profit(t,T)$ with high accuracy: in both cases the coefficient of determination $R^2$ was very close to 1.

Comparisons of metamodels for $Purity(t,T)$ and $Profit(t,T)$, $i$=1, 2 (26) with the full model solutions show a very good agreement between surfaces (Fig. 5, a, c). DS curves (red lines) defined by $DS_{(\boldsymbol{u};Purity(\boldsymbol{u}))} = (\boldsymbol{u}, Purity = P^*)$, $DS_{(\boldsymbol{u};Profit(\boldsymbol{u}))} = (\boldsymbol{u}, Profit = Pr^*)$ also agree very well (Fig. 5, b, d).

Metamodels for CQA's and corresponding DS functions are given explicitly in analytical forms (26). It allows to apply the theory of R-functions for identification of a joint DS from equations for individual DS's. Following the framework presented in Section 5.2, we define implicit functions

$$\phi_1(T,t) = -P^* + a_{00}^{[1]} + a_{10}^{[1]}T + a_{20}^{[1]}T^2 + a_{01}^{[1]}t + a_{11}^{[1]}t\,T,$$

$$\phi_2(T,t) = -Pr^* + a_{00}^{[2]} + a_{10}^{[2]}T + a_{20}^{[2]}T^2 + a_{01}^{[2]}t + a_{11}^{[2]}t\,T. \tag{27}$$

Their values for the considered model are given by Eqs. (28), (29). Fig. 6 shows the DS curves for $Purity(t,T)$ and $Profit(t,T)$ defined by (28), (29).

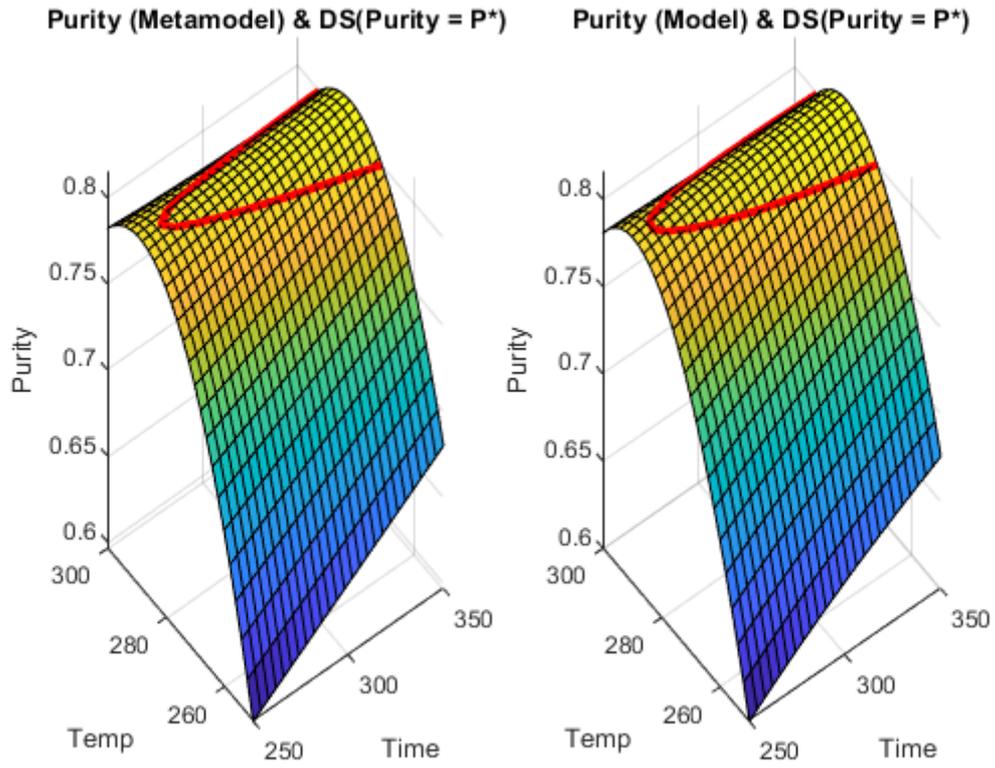

(a)

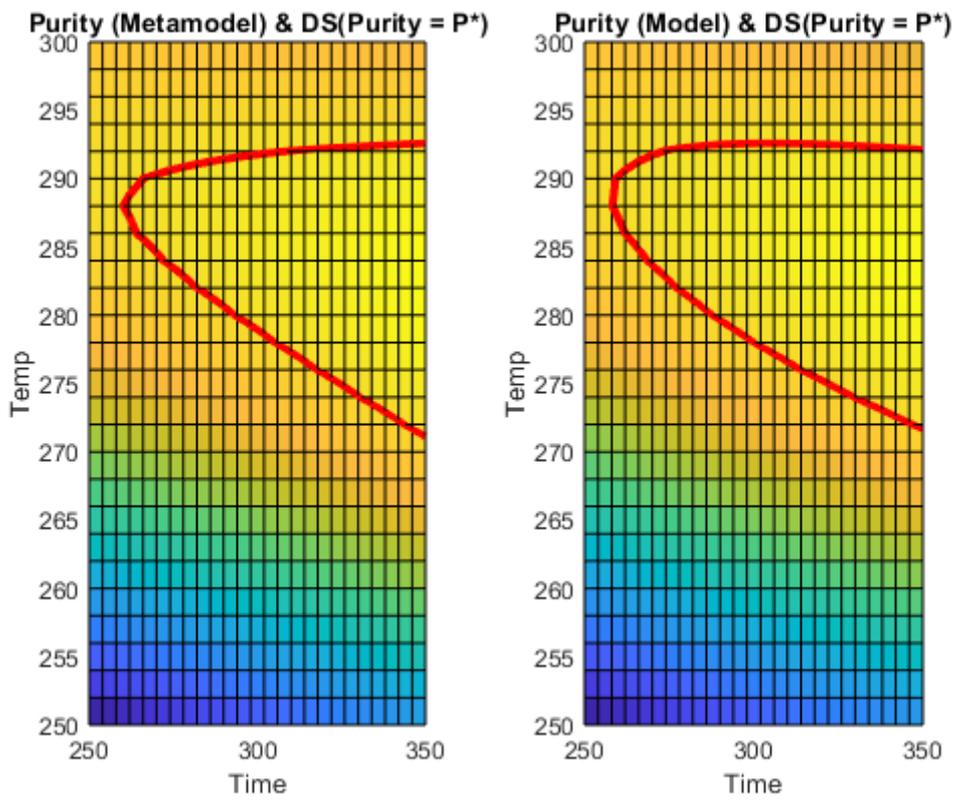

(b)

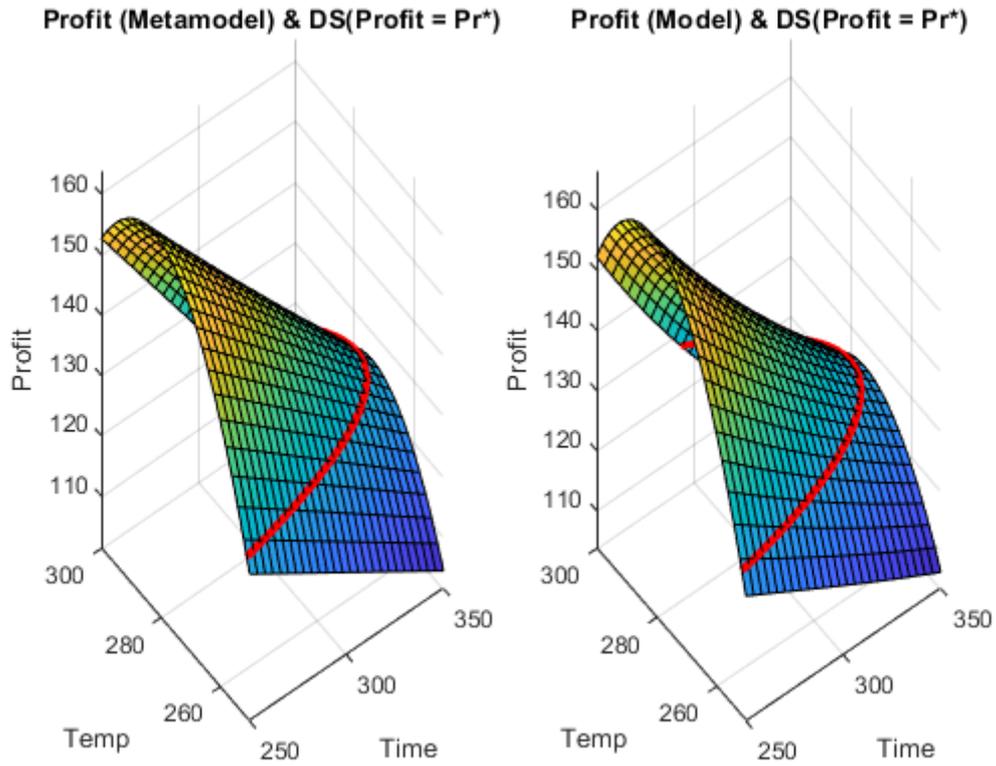

(c)

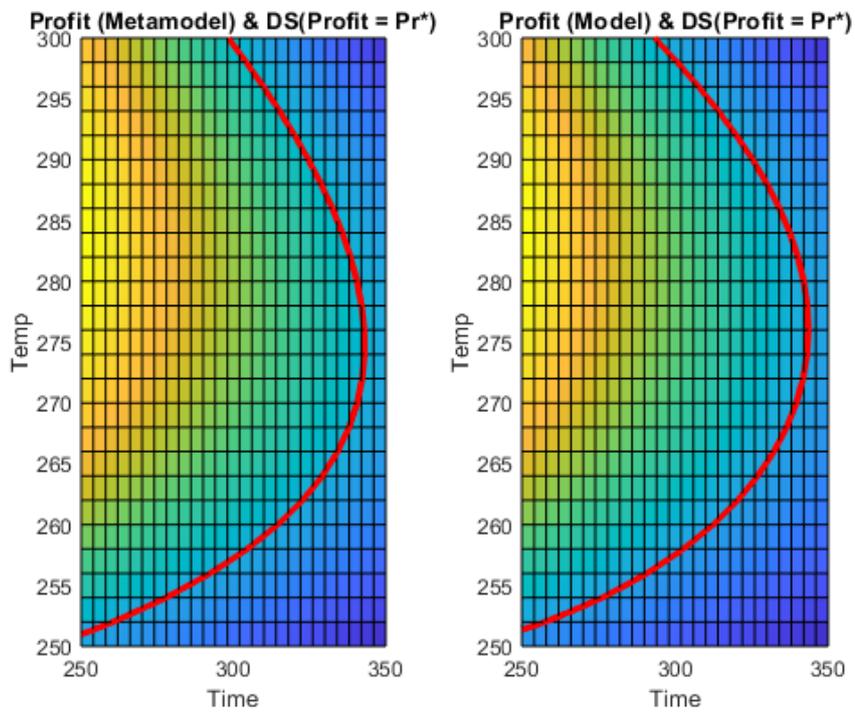

(d)

Fig. 5. $Purity(t,T)$. Metamodel (20) (left) and full model based solutions (right) 3D plot (a); $Purity(t,T)$. Metamodel (20) (left) and full model based solutions (right). 2D projections (b); $Profit\ (t,T)$. Metamodel (20) (left) and full model based solutions (right). 3D plot (c);

*Profit* $(t, T)$. Metamodel (20) (left) and full model based solutions (right). 2D projections (d); $DS_{(\mathbf{u};Purity(\mathbf{u}))} = (\mathbf{u}, Purity = P^*)$, $DS_{(\mathbf{u};Profit(\mathbf{u}))} = (\mathbf{u}, Profit = Pr^*)$ curves in red;

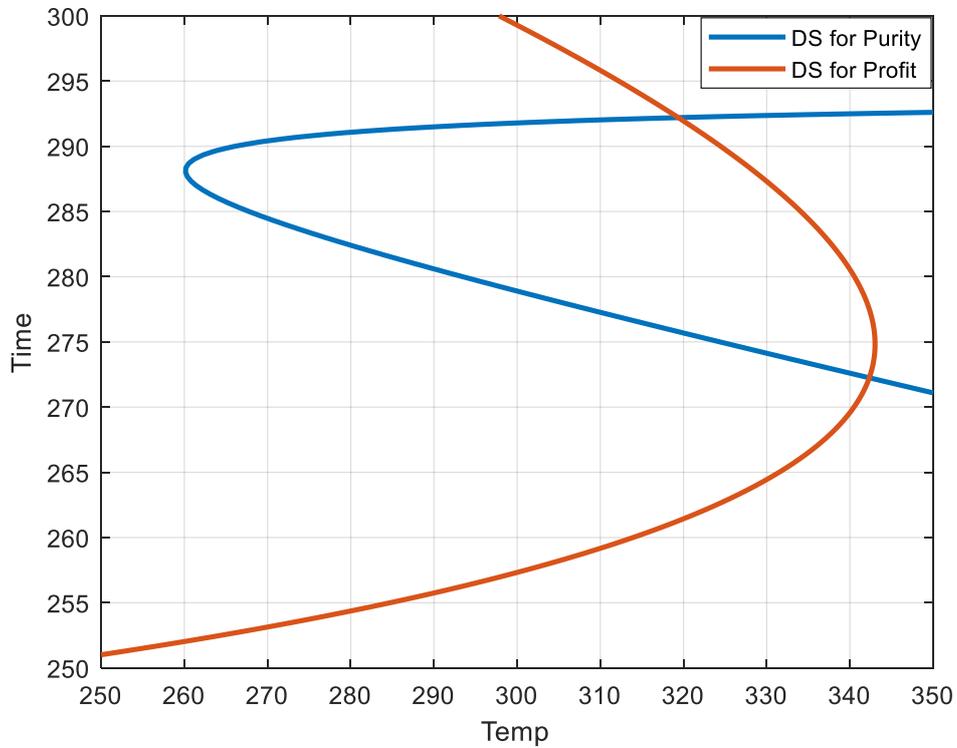

Fig. 6. DS curves for *Purity* (blue) and *Profit* (red).

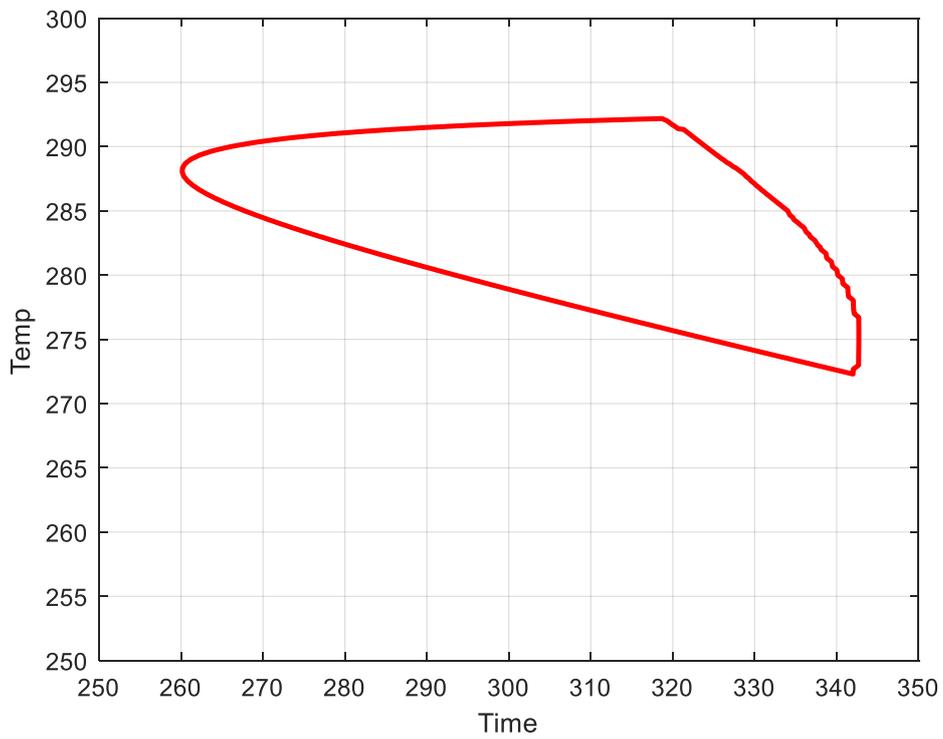

(a)

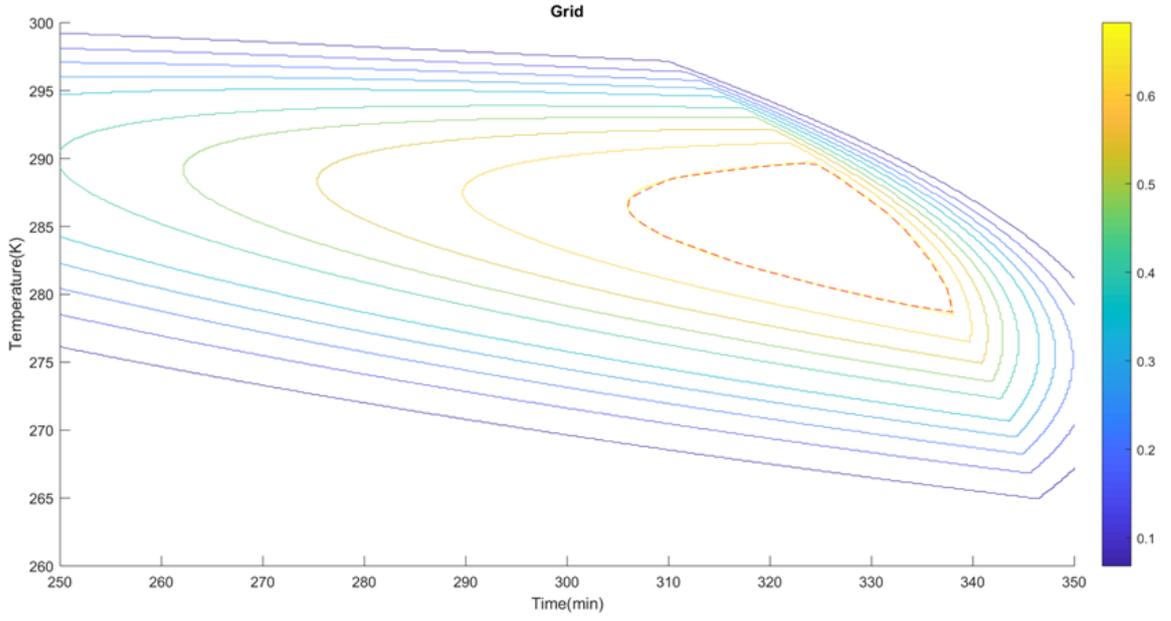

(b)

Fig. 7. Analytically defined (Eq. (30)) joint DS (a); Contour plots of the probability map obtained using an "exhaustive sampling" method (Kucherenko et al. 2019) (b).

$$\phi_1(T,t) = \frac{1492155457490699 * T}{18014398509481984} + \frac{6353038647553177 * t}{1152921504606846976} \\ - \frac{5532151810052315 * T * t}{295147905179352825856} - \\ \frac{2496177142651283 * T^2}{18446744073709551616} - \frac{8909117195501237}{703687441776640}. \quad (28)$$

$$\phi_2(T,t) = \frac{6181541458134125 * T}{281474976710656} + \frac{1444357484257337 * t}{1125899906842624} \\ - \frac{6903876884140963 * T * t}{1152921504606846976} \\ \frac{(5220857310952963 * T^2)}{144115188075855872} - \frac{6981627641092835}{2199023255552}. \quad (29)$$

Applying R-function $f_\wedge(\phi_1, \phi_2)$ we obtain the joint DS in an analytical form:.

$$f_\wedge(T,t) = \frac{397110808778074699 * T}{36028797018963968} + \frac{1485375102527066265 * t}{2305843009213693952} - \frac{17729246341501388843\, T*T}{590295810358705651712}$$

$$\frac{-((8809301652650171065 * T * t - 18850163522969086080 * t}{}$$
$$\frac{-32286842704944603217920 * T + 53261884692746238480 * T^2 +}{}$$
$$\frac{46666072257026189437173 76)^2)^{(\frac{1}{2})}}{2951479051793528258560}$$

$$-\frac{(670765912944630547 * T^2)}{36893488147419103232} - \frac{2243029962345208437}{1407374883553280}.$$

(30)

The joint DS defined by Eq. (30) is shown in Fig. 7 (a). A comparison between the DS obtained using R-function for a deterministic model (Fig. 7 (a)) with that obtained using the acceptance – rejection method sampling approach (Kucherenko et al. 2019) for a probabilistic model (Fig. 7(b)) shows a good agreement between the two methods. The R-DS identifier method required only $N$=64 model runs to achieve the same accuracy as the acceptance – rejection method, which required $N$ = 1640 model runs, while the "exhaustive sampling" method required $2.6 \cdot 10^5$ model runs. Here the number of model runs $N$ accounts only for model runs needed to produce a deterministic part of the probability maps, excluding the probabilistic part. In the acceptance – rejection and "exhaustive sampling" methods the solution for the DS is given a set of "accepted points". It means that the problem of the DS boundary description in a practical, usable form remains unsolved, while the R-DS identifier method gives an explicit analytical expression for the DS boundary.

The proposed method is relatively easy to implement using symbolic computation tools available in software platforms like MATLAB (The MathWorks, 2023), Maple, and Mathematica. They enable the derivation of implicit functions $\phi_i(\boldsymbol{u})$ for individual DS's and the joint DS $f_\wedge(\boldsymbol{u})$. In this work we used MATLAB for all symbolic computations and plotting.

## 6. Conclusions

Rvachev's R-functions are incredibly useful for representing geometric shapes and operations. Their ability to combine implicit functions into complex shapes through algebraic expressions makes them invaluable in fields ranging from computer graphics to engineering analysis. We showed that they also enable the solution of geometric problems such as identification of the Design Space through algebraic manipulation.

Pharmaceutical manufacturing model systems are usually very complex, making the task of defining the Design Space computationally intractable with existing methods. We developed a new methodology for efficient Design Space identification by using R-functions to describe the Design Space in a form of analytical expressions. This approach can effectively convey high-dimensional Design Space shapes and incorporate this information into the formal documentation required for Design Space approval by a government agency.

## Acknowledgements

S.K and N.S. acknowledge the financial support of Eli Lilly and Company and the EPSRC Grant EP/T005556/1.

## Appendix
### Test cases in 3D. Planes

Consider three equations defining infinite 3D slabs of thickness $2a$, $2b$ and $2c$ surrounding each of the coordinate planes $(y, z)$, $(x, z)$ and $(x, y)$, respectively (Fig. A.1 (a)):

$$f_1(x, y, z) = a^2 - x^2, f_2(x, y, z) = b^2 - y^2, f_3(x, y, z) = c^2 - z^2. \quad (A.1)$$

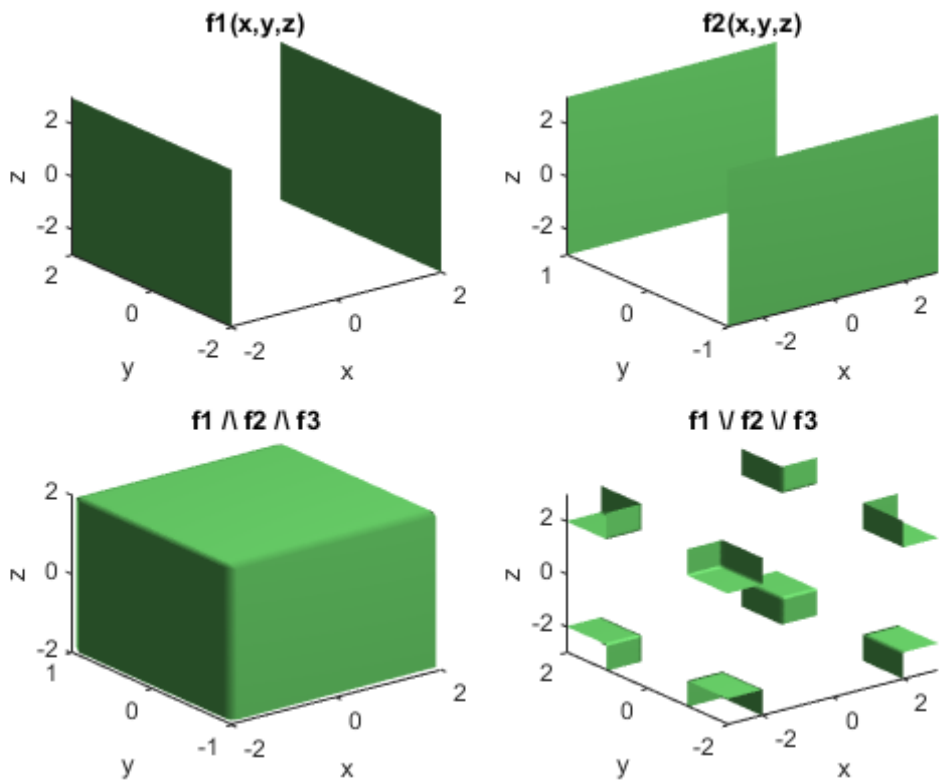

Fig. A.1. $f_1(x, y, z) = 0$ (a); $f_2(x, y, z) = 0$ (b); **R**-conjunction of all three functions in (A.1) (c); **R**-disjunction of all three functions in (A.1) (d).

Following the procedure outlined in the previous Section, we build a **R**-conjunction $f_\wedge(x, y, z) = f_\wedge(\phi_1, \phi_2, \phi_3)$ and a **R**-disjunction $f_\vee(x, y, z) = f_\vee(\phi_1, \phi_2, \phi_3)$ of the three functions (Fig. A.1) assuming $a = 2$, $b = 1$ and $c = 2$:

$$f_\wedge(x,y,z) = 9 - x^2 - y^2 - z^2 - g(x,y) - \left(\big(g(x,y) + h(x,y)\big)^2 - (z^2 - 4)\big(g(x,y) + h(x,y)\big) + (z^2 - 4)^2\right)^{\frac{1}{2}},$$

$$f_\vee(x,y,z) = 9 - x^2 - y^2 - z^2 + g(x,y) + \left(\big(g(x,y) - h(x,y)\big)^2 + (z^2 - 4)\big(g(x,y) - h(x,y)\big) + (z^2 - 4)^2\right)^{\frac{1}{2}},$$

where $g(x,y) = \left((x^2 - 4)^2 + (y^2 - 1)^2 - (x^2 - 4)(y^2 - 1)\right)^{\frac{1}{2}}$ and $h(x,y) = x^2 + y^2 - 5$.

The domain described by the inequality $f_\wedge(x,y,z) \geq 0$ is a 3D cuboid (Fig. A.1 (c)) and its whole boundary is implicitly defined by a single equation $f_\wedge(x,y,z) = 0$.

The R-disjunction of the three slabs represents the complement of the cuboid in Fig. A.1 (c), that is eight octants separated by finite gaps aligned with the coordinate planes (Fig. A.1 (d)).

### Test cases in 3D. Paraboloids and cylinders

Consider the following inequalities that define:

the area under the paraboloid $z = 0.6(1 - x^2 - y^2)$:

$$f_1(x,y,z) = -z + 0.6(1 - x^2 - y^2) \geq 0 \tag{A.2}$$

the area above the paraboloid $z = -0.6(1 - x^2 - y^2)$:

$$f_2(x,y,z) = z + 0.6(1 - x^2 - y^2) \geq 0 \tag{A.3}$$

the area inside of the cylinder of radius 0.5 centered around the z-axis:

$$f_3(x,y,z) = 0.5^2 - x^2 - y^2 \geq 0 \tag{A.4}$$

the area inside of the cylinder of radius 0.3 centered around the z-axis:

$$f_4(x,y,z) = 0.3^2 - x^2 - y^2 \geq 0 \tag{A.5}$$

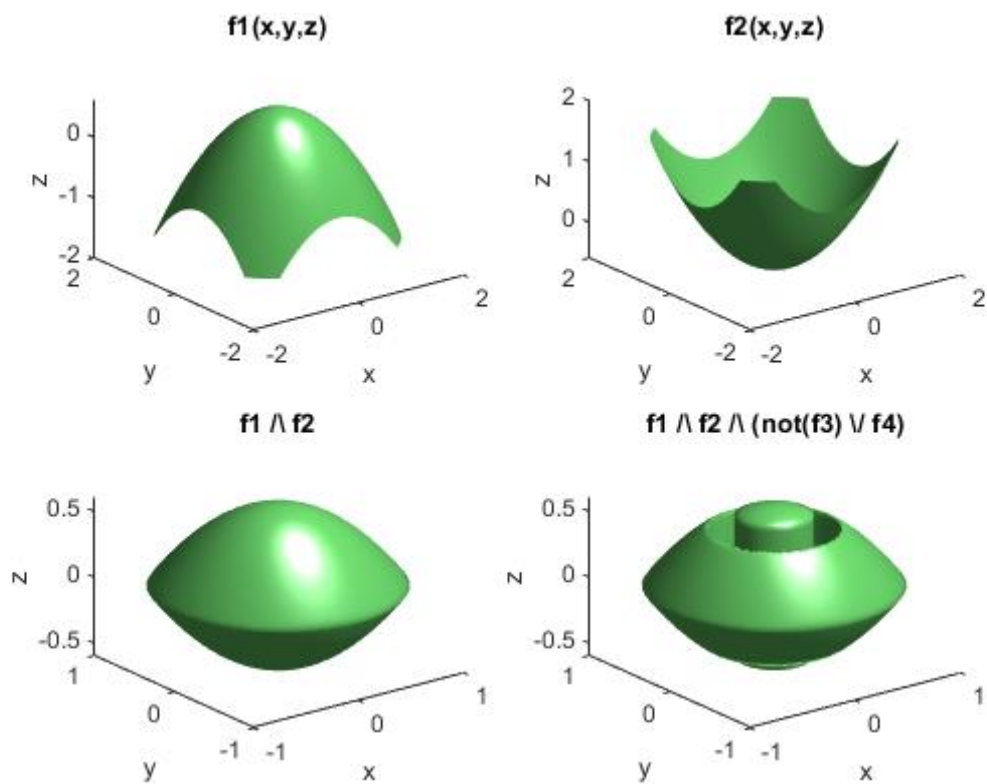

Fig. A.2. $f_1(x, y, z) = 0$ (a); $f_2(x, y, z) = 0$ (b); **R**-conjunction of $f_1(x, y, z)$ and $f_2(x, y, z)$ (c); **R**-conjunction of $f_1(x, y, z)$, $f_2(x, y, z)$ and $f_1 \wedge_1 f_2 \wedge_1 (\overline{f_3} \vee_1 f_4)$ (d).

Using these geometrical primitives, we can construct domains either using simple operations as above, such as an R-conjunction of $f_1(x, y, z)$ and $f_2(x, y, z)$ (Fig. A.2(c)), as well as their combinations to describe complex non-simply connected domains like the one shown in Fig. A.2 (d). The latter is described by the R-function $f_1 \wedge_1 f_2 \wedge_1 (\overline{f_3} \vee_1 f_4) \geq 0$ and represents the domain shown in Fig. A.2 (c) with an annular cylindrical cut-out.